\documentclass[11pt]{article}

\usepackage{amsmath}
\usepackage{url}
\usepackage{amsthm}

\usepackage[margin=1.5in]{geometry}

\usepackage{float}

\usepackage{graphicx}
\usepackage{psfrag}

\usepackage[format=hang,indention=0pt,figurename=Fig.]{caption}

\usepackage{amsfonts}

\newcommand{\nR}{\mathbb R}
\newcommand{\nRc}{{\mathbb R}\cup\{\infty\}}

\newtheorem*{theorem*}{Theorem}
\newtheorem*{corollary*}{Corollary}
\begin{document}

\title{A generalization of Routh's triangle theorem}

\author{\'Arp\'ad B\'enyi and Branko \'Curgus}


\maketitle

\begin{abstract}
We prove a generalization of the well known Routh's triangle theorem. As a consequence, we get a unification of the theorems of Ceva and Menelaus. A connection to Feynman's triangle is also given.
\end{abstract}

\section{Introduction.}

On page 82 in \cite{Routh}, just after his now famous triangle theorem, Routh writes:

\begin{center}
\begin{minipage}[c]{.85\textwidth}
``The author has not met with these expressions for the area of two
triangles which often occur. He has therefore placed them here in
order that the argument in the text may be more easily understood.''
\end{minipage}
\end{center}
In this note we revisit ``these expressions'' of Routh and show that they can be thought of as special cases of only one expression.  To our surprise, this unifying expression relating the areas of two triangles seems to be missing from the literature.  Thus, paraphrasing Routh, we deemed it useful to place it here. We hope that the reader will find our generalization of Routh's triangle theorem not only natural but also intrinsically beautiful in its symmetry.

Let $ABC$ be a triangle determined by three non-collinear points  $A, B, C$ in the Euclidean plane. A line joining a vertex to a point on the line containing the opposite side is called a {\em cevian}. Given a point $D$ on the line $BC$, the line $AD$ is a cevian through the vertex $A$. If $D \neq C$, then this cevian is uniquely determined by $x \in \nR\!\setminus\!\{-1\}$ such that $\overrightarrow{BD} =  x \overrightarrow{DC}$.  Conversely, for $x \in \nR\!\setminus\!\{-1\}$, the equality $\overrightarrow{BD} =  x \overrightarrow{DC}$ determines uniquely a point $D \neq C$ on the line $BC$.  In this note, for such a point $D$, we write $D = A_x$. It follows directly from the definition that $A_0 = B$. We adopt the notation $A_{\infty} = C$. The line through $A$ which is parallel to $BC$ will also be considered as a cevian through the vertex $A$.  We denote it by $AA_{-1}$, thinking of $A_{-1}$ as the ``point at infinity'' on the line $BC$.  In this way, we establish a bijection between the points of the line $BC$ and the set $\nRc$. This bijection associates positive numbers to the points between $B$ and $C$, it associates numbers in $(-1,0)$ to the points between the point at infinity and $B$, and it associates numbers less then $-1$ to the points between the point at infinity and $C$. We use an analogous notation for the cevians through vertices $B$ and $C$. We will denote by  $B_y$  the unique point on the line $CA$ such that $\overrightarrow{CB_y} = y\overrightarrow{B_yA}$ and by $C_z$ the unique point on the line $AB$ such that $\overrightarrow{AC_z} = z \, \overrightarrow{C_zB}$. By definition, $B_0 = C, B_\infty = A, C_0 = A$ and $C_\infty = B$.

\section{Two expressions of Routh.}

The first of the two expressions that Routh refers to in the displayed quote gives the ratio between the area $\Delta'$ of the triangle $A_xB_yC_z$ and the area $\Delta$ of the triangle $ABC$, see \cite[p.~82]{Routh}:
\begin{equation} \label{eqDc}
\frac{\Delta'}{\Delta} = \frac{xyz + 1}{(1+x)(1+y)(1+z)}.
\end{equation}
The natural domain of \eqref{eqDc} is $(\nR\setminus\{-1\})^3$. In fact, it is an exercise in multivariable limits that \eqref{eqDc} can be continuously extended to $\bigl((\nRc) \setminus\{-1\}\bigr)^3$. In this way, \eqref{eqDc} holds whenever $A_x, B_y, C_z$, with $x,y,z\in \nRc$, are ``finite points''.

In \eqref{eqDc} and all subsequent formulas involving areas we  consider signed areas of triangles. A positive area corresponds to a positively oriented triangle; that is a triangle in which the increasing alphabetic order of vertices proceeds counterclockwise. A negative area corresponds to a triangle of opposite orientation.

If $x=y=z=1$, then the triangle $A_xB_yC_z$ is commonly known as the {\em medial triangle} of $ABC$. In analogy with this terminology, we will refer to the general triangle $A_xB_yC_z$ as a {\em cevial triangle} of $ABC$. It follows from \eqref{eqDc} that the points $A_x, B_y, C_z$ are collinear if and only if $xyz = -1$. This is the well known theorem of Menelaus, see \cite[p.~220]{Coxeter}.

The second of the two expressions that Routh refers to in the quote is another equality, nowadays known as Routh's theorem. Whenever the cevians $AA_x$ and $BB_y$ intersect at a single point, we denote by $P$ their intersection point. Similarly, we denote by $Q$ the intersection point of $BB_y$ and $CC_z$ and by $R$ the intersection point of $CC_z$ and $AA_x$.  Then, Routh's theorem gives the ratio between the area $\Delta''$ of the triangle $PQR$ and the area $\Delta$  of the triangle $ABC$, see again \cite[p.~82]{Routh}:
\begin{equation} \label{eqDr}
\frac{\Delta''}{\Delta} = \frac{(xyz - 1)^2}{(1 + x + x y)(1 + y + y z)(1 + z + z x)}.
\end{equation}
The natural domain of \eqref{eqDr} is $\nR^3\!\setminus\!S$, where $S$ is the set of all triples $(x,y,z) \in \nR^3$ such that $(1+x+xy)(1+y+yz)(1+z+zx) = 0$.  As before, \eqref{eqDr} can be continuously extended to $(\nRc)^3\!\setminus\!\overline{S}$, where $\overline{S}$ is the closure of $S$ in $(\nRc)^3$.  To understand the set $\overline{S}$, notice that the only solutions of $1+x+xy = 0$ involving $\infty$ are $(x,y) = (\infty,-1)$ and $(x,y) = (0, \infty)$.  This extended domain of \eqref{eqDr} coincides with the set of triples $(x,y,z) \in (\nRc)^3$ for which each of the pairs of cevians  $(AA_x,BB_y)$, $(BB_y,CC_z)$ and $(CC_z,AA_x)$ intersects at exactly one point. More specifically, with $x,y \in \nRc$, the cevians $AA_x$ and $BB_{y}$ are parallel if and only if $1+x+xy =0$ and $(x,y) \neq (0, \infty)$; this follows from the calculations preceding \eqref{eqPQR} below. Note also that the cevians $AA_0$ and $BB_{\infty}$ coincide.

The triangle $PQR$ will be called a {\em Routh's triangle} of $ABC$. Notice that the points $P, Q$ and $R$ are collinear if and only if they coincide. Therefore, \eqref{eqDr} implies that the lines $AA_x, BB_y$ and $CC_z$ are concurrent if and only if $xyz = 1$. This statement is known as Ceva's theorem, \cite[p.~220]{Coxeter}.

\section{A unification.}

\begin{figure}[H]

\begin{center}
\psfrag{C}[][]{\begin{picture}(0,0)
            \put(-5,8){\makebox(0,0)[l]{$C$}}
                        \end{picture}}
\psfrag{Bv}[][]{\begin{picture}(0,0)
            \put(-16,5){\makebox(0,0)[l]{$B_v$}}
                        \end{picture}}
\psfrag{Ax}[][]{\begin{picture}(0,0)
            \put(0,6){\makebox(0,0)[l]{$A_x$}}
                        \end{picture}}
\psfrag{rA}[][]{\begin{picture}(0,0)
            \put(-6,9){\makebox(0,0)[l]{$P$}}
                             \end{picture}}
\psfrag{By}[][]{\begin{picture}(0,0)
            \put(-16,5){\makebox(0,0)[l]{$B_y$}}
                        \end{picture}}
\psfrag{rB}[][]{\begin{picture}(0,0)
            \put(-10,-4){\makebox(0,0)[l]{$Q$}}
                        \end{picture}}
\psfrag{Au}[][]{\begin{picture}(0,0)
            \put(0,6){\makebox(0,0)[l]{$A_u$}}
                        \end{picture}}
\psfrag{rC}[][]{\begin{picture}(0,0)
            \put(-10,-8){\makebox(0,0)[l]{$R$}}
                        \end{picture}}
\psfrag{A}[][]{\begin{picture}(0,0)
            \put(-8,-8){\makebox(0,0)[l]{$A$}}
                        \end{picture}}
\psfrag{Cw}[][]{\begin{picture}(0,0)
            \put(-4,-9){\makebox(0,0)[l]{$C_w$}}
                        \end{picture}}
\psfrag{Cz}[][]{\begin{picture}(0,0)
            \put(-4,-9){\makebox(0,0)[l]{$C_z$}}
                        \end{picture}}
\psfrag{B}[][]{\begin{picture}(0,0)
            \put(-2,-8){\makebox(0,0)[l]{$B$}}
                        \end{picture}}

\setlength{\abovecaptionskip}{-5pt}%
\setlength{\belowcaptionskip}{-10pt}%

 \resizebox{!}{!}{\includegraphics{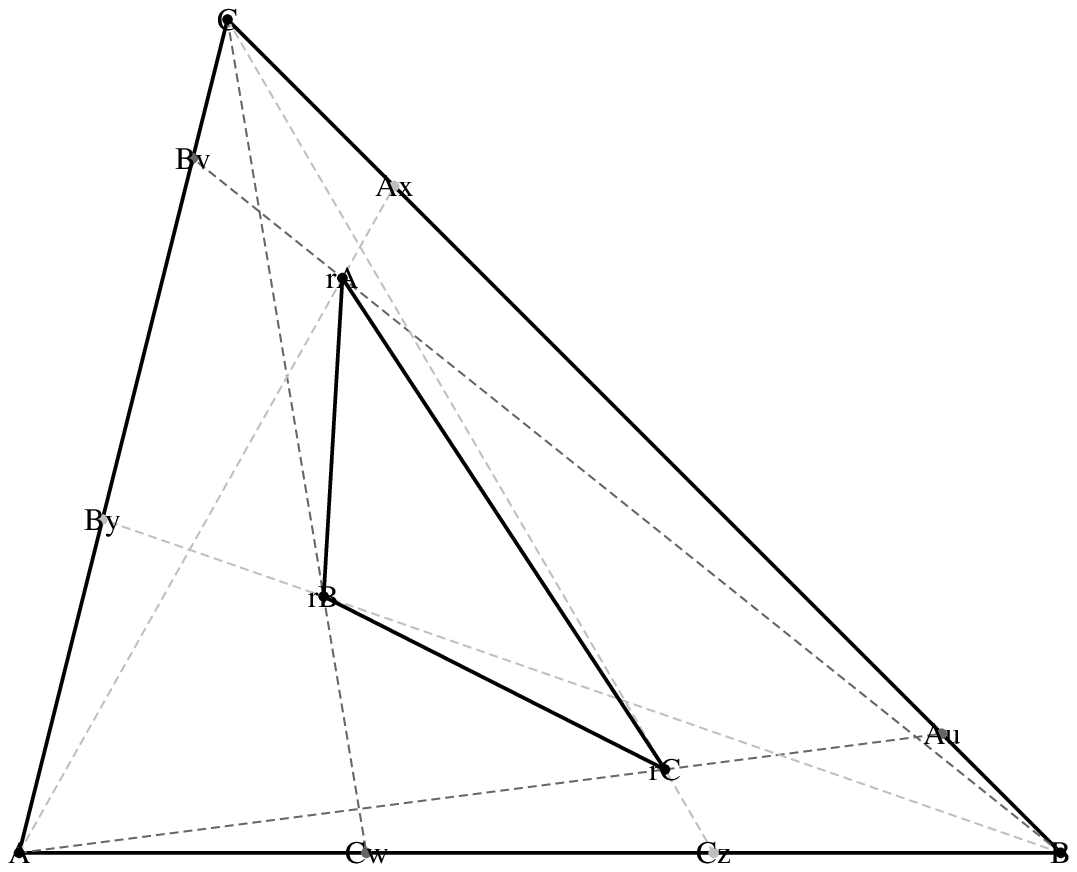}}
  \caption{A generalized Routh's triangle} \label{f1}
\end{center}
\end{figure}

An extension of the construction by Nakamura and Oguiso
\cite[\S3]{Nakamura} allows for formulas \eqref{eqDc} and \eqref{eqDr} to be unified. For $u, v, w, x, y, z \in \nRc$, consider the following six cevians, two from each vertex: $AA_u, AA_x, BB_v, BB_y, CC_w, CC_z$. Assuming that each of the pairs of cevians $(AA_x,BB_v)$, $(BB_y,CC_w)$ and $(CC_z,AA_u)$ intersects at exactly one point, we define the {\em generalized Routh's triangle} $PQR$ of the given triangle $ABC$, see Figure~\ref{f1}, by setting
\begin{equation} \label{eqgKLM}
 \{P\} =  AA_x \cap BB_v, \qquad  \{Q\} = BB_y \cap CC_w, \qquad  \{R\} = CC_z \cap AA_u.
\end{equation}

\vspace*{-35pt}

\begin{figure}[H]
\begin{minipage}{0.49\textwidth}
\psfrag{C}[][]{\begin{picture}(0,0)
            \put(-5,8){\makebox(0,0)[l]{$C$}}
                        \end{picture}}
\psfrag{Bv}[][]{\begin{picture}(0,0)
            \put(-17,0){\makebox(0,0)[l]{$B_v$}}
                        \end{picture}}
\psfrag{Ax}[][]{\begin{picture}(0,0)
            \put(3,4){\makebox(0,0)[l]{$A_x$}}
                        \end{picture}}
\psfrag{rA}[][]{\begin{picture}(0,0)
            \put(-6,9){\makebox(0,0)[l]{$P$}}
                             \end{picture}}
\psfrag{By}[][]{\begin{picture}(0,0)
            \put(-15,7){\makebox(0,0)[l]{$B_y$}}
                        \end{picture}}
\psfrag{rB}[][]{\begin{picture}(0,0)
            \put(0,6){\makebox(0,0)[l]{$Q$}}
                        \end{picture}}
\psfrag{Au}[][]{\begin{picture}(0,0)
            \put(-4,10){\makebox(0,0)[l]{$A_u$}}
                        \end{picture}}
\psfrag{rC}[][]{\begin{picture}(0,0)
            \put(-12,0){\makebox(0,0)[l]{$R$}}
                        \end{picture}}
\psfrag{A}[][]{\begin{picture}(0,0)
            \put(-8,-8){\makebox(0,0)[l]{$A$}}
                        \end{picture}}
\psfrag{Cw}[][]{\begin{picture}(0,0)
            \put(1,-9){\makebox(0,0)[l]{$C_w$}}
                        \end{picture}}
\psfrag{Cz}[][]{\begin{picture}(0,0)
            \put(-9,-9){\makebox(0,0)[l]{$C_z$}}
                        \end{picture}}
\psfrag{B}[][]{\begin{picture}(0,0)
            \put(-2,-8){\makebox(0,0)[l]{$B$}}
                        \end{picture}}

\setlength{\abovecaptionskip}{2pt}%
\setlength{\belowcaptionskip}{-6pt}%

\resizebox{\textwidth}{!}{\includegraphics{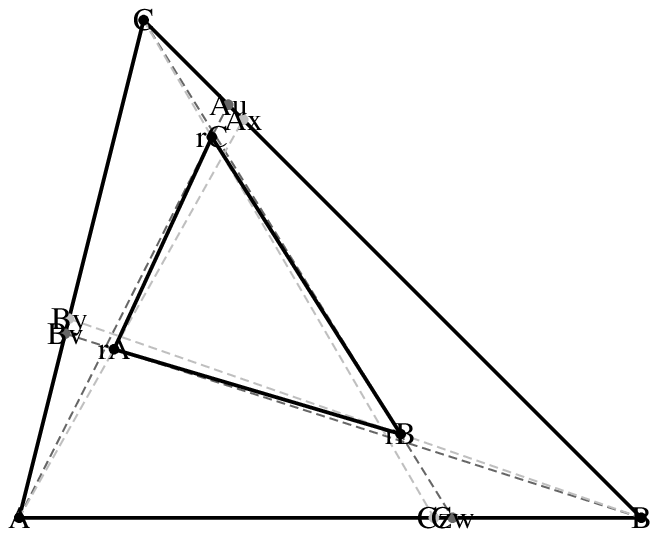}}
  \caption{Almost Routh's $PQR$} \label{f2}
\end{minipage}
\begin{minipage}{0.49\textwidth}
\psfrag{C}[][]{\begin{picture}(0,0)
            \put(-5,8){\makebox(0,0)[l]{$C$}}
                        \end{picture}}
\psfrag{Bv}[][]{\begin{picture}(0,0)
            \put(-16,5){\makebox(0,0)[l]{$B_v$}}
                        \end{picture}}
\psfrag{Ax}[][]{\begin{picture}(0,0)
            \put(0,6){\makebox(0,0)[l]{$A_x$}}
                        \end{picture}}
\psfrag{rA}[][]{\begin{picture}(0,0)
            \put(-13,1){\makebox(0,0)[l]{$P$}}
                             \end{picture}}
\psfrag{By}[][]{\begin{picture}(0,0)
            \put(-16,5){\makebox(0,0)[l]{$B_y$}}
                        \end{picture}}
\psfrag{rB}[][]{\begin{picture}(0,0)
            \put(6,2){\makebox(0,0)[l]{$Q$}}
                        \end{picture}}
\psfrag{Au}[][]{\begin{picture}(0,0)
            \put(0,6){\makebox(0,0)[l]{$A_u$}}
                        \end{picture}}
\psfrag{rC}[][]{\begin{picture}(0,0)
            \put(0,7){\makebox(0,0)[l]{$R$}}
                        \end{picture}}
\psfrag{A}[][]{\begin{picture}(0,0)
            \put(-8,-8){\makebox(0,0)[l]{$A$}}
                        \end{picture}}
\psfrag{Cw}[][]{\begin{picture}(0,0)
            \put(-4,-9){\makebox(0,0)[l]{$C_w$}}
                        \end{picture}}
\psfrag{Cz}[][]{\begin{picture}(0,0)
            \put(-4,-9){\makebox(0,0)[l]{$C_z$}}
                        \end{picture}}
\psfrag{B}[][]{\begin{picture}(0,0)
            \put(-2,-8){\makebox(0,0)[l]{$B$}}
                        \end{picture}}

\setlength{\abovecaptionskip}{-2pt}%
\setlength{\belowcaptionskip}{-6pt}%

\resizebox{\textwidth}{!}{\includegraphics{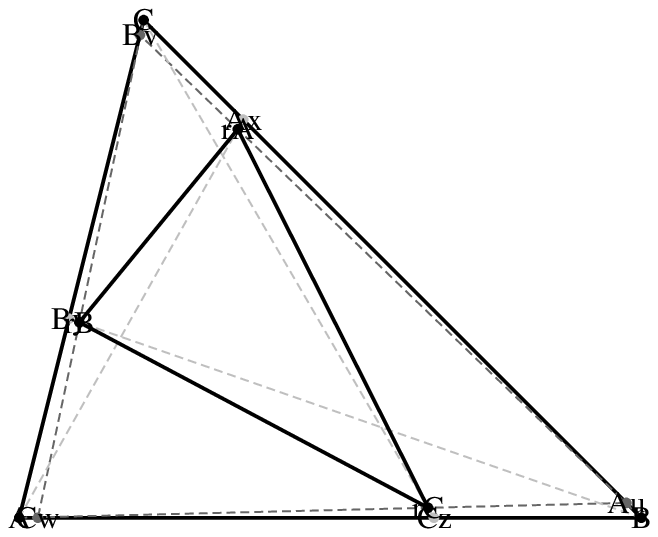}}
  \caption{Almost a cevial $PQR$} \label{f3}
\end{minipage}

\end{figure}

Note that the discussion in the paragraph following \eqref{eqDr} implies that the points $P, Q$ and $R$ in \eqref{eqgKLM} are well defined if and only if $x,y,z,u,v,w \in \nRc$ satisfy
\begin{equation} \label{eqGRc}
(1+x+xv)(1+y+yw)(1+z+zu) \neq 0.
\end{equation}
Choosing $u=x, v=y,$ and $w=z$ the triangle $PQR$ defined in \eqref{eqgKLM} becomes a Routh's triangle, see Figure~\ref{f2}. On the other hand, choosing $u=v=w=0$ we have $AA_0 = AB, BB_0 = BC, CC_0 = CA$, and  \eqref{eqgKLM} yields $P = A_x, Q = B_y, R = C_z$. In this case, the triangle $PQR$ is the cevial triangle $A_xB_yC_z$, see Figure~\ref{f3}.

\vspace*{-15pt}

\begin{theorem*}
With the points $P, Q, R$ defined in \eqref{eqgKLM}, the ratio between the  area $\Delta_1$ of the triangle $PQR$ and the area $\Delta$ of the triangle $ABC$ is given by
\begin{equation} \label{eqD}
\frac{\Delta_1}{\Delta} = \frac{1-x y w - x v z - u y z + x y z + xyzuvw}%
{(1+x+xv)(1+y+yw)(1+z+zu)}.
\end{equation}
\end{theorem*}

The natural domain of \eqref{eqD} is the set of all $(x,y,z,u,v,w) \in \nR^6$ which satisfy \eqref{eqGRc}. However, it is again an exercise in multivariable limits, this time with six variables, to check that formula \eqref{eqD} extends by continuity to all $(x,y,z,u,v,w) \in (\nRc)^6$ which satisfy \eqref{eqGRc}.  Beautifully, with $u = x, v = y, w = z$, \eqref{eqD} simplifies to \eqref{eqDr} and with $u = v = w = 0$, it simplifies to \eqref{eqDc}. Also, as $x,y,z \to \infty$, formula \eqref{eqD} becomes \eqref{eqDc}, with $x$ in \eqref{eqDc} substituted by $u$, $y$ by $v$, and $z$ by $w$.

It is quite possible that any of the many proofs of Routh's theorem (see, for example, \cite[Section~13.7]{Coxeter}, \cite{Kline}, \cite{Niven}, to mention a few) can be modified to prove this generalization.  We will prove it by using two  standard undergraduate tools: linear algebra and analytic geometry. First, we observe that the ratio of the areas remains unchanged under affine transformations of $\nR^2$. Therefore, instead of an arbitrary triangle $ABC$, we can consider the triangle with the vertices
\begin{equation*}
A = (0,0), \quad B = (1,0), \quad \text{and} \quad \quad C = (0,1).
\end{equation*}
Next, we find the coordinates of the point $P$ for these special vertices $A,B,C$. Let $\xi$ denote its first coordinate and calculate $A_x = \bigl(1/(1+x), x/(1+x)\bigr)$ and $B_v = \bigl(0,1/(1+v)\bigr)$. Now, intersecting the lines $AA_x$ and $BB_v$ we get $\xi x = (1-\xi)/(1+v)$. Solving for $\xi$ gives the first coordinate of $P$, while the second coordinate is $\xi x$. Similarly, we find $Q$ and $R$; we have
\begin{equation} \label{eqPQR}
P  = \bigl(\tfrac{1}{1+x+xv},\tfrac{x}{1+x+xv}\bigr),  \ \
Q   = \bigl(\tfrac{yw}{1+y+yw},\tfrac{1}{1+y+yw}\bigr), \ \
R   = \bigl(\tfrac{z}{1+z+zu},\tfrac{zu}{1+z+zu}\bigr).
\end{equation}
Let $\Box$ stand for the left-hand side of \eqref{eqGRc}. Since we assume that each pair of cevians $\bigl(AA_x,BB_v\bigr)$, $\bigl(BB_y,CC_w\bigr)$ and $\bigl(CC_z,AA_u\bigr)$ intersects at exactly one point, we have $\Box \neq 0$. As the area of the triangle $ABC$ is $1/2$, the ratio in the theorem is given by the determinant
 \allowdisplaybreaks{
\begin{align*}
\left| \begin{array}{ccc}
 \frac{1}{1+x+xv} & \frac{yw}{1+y+yw} & \frac{z}{1+z+zu} \\[6pt]
 \frac{x}{1+x+xv} & \frac{1}{1+y+yw} &  \frac{zu}{1+z+zu} \\[6pt]
1 & 1 & 1
 \end{array} \right|
 & = \frac{1}{\Box}
 \left| \begin{array}{ccc}
 1 & yw & z \\[6pt]
 x & 1  & zu \\[6pt]
1+x+xv & 1+y+yw & 1+z+zu
 \end{array} \right| \\[8pt]
 & \hspace*{-1.05in} = \frac{1}{\Box}
  \left| \begin{array}{ccc}
 1 & yw & z \\[6pt]
 x & 1  & zu \\[6pt]
1  & 1  & 1
 \end{array} \right| + \frac{1}{\Box}
 \left| \begin{array}{ccc}
 1 & yw & z \\[6pt]
 x & 1  & zu \\[6pt]
 x  &  y  &  z
 \end{array} \right| + \frac{1}{\Box}
 \left| \begin{array}{ccc}
 1 & yw & z \\[6pt]
 x & 1  & zu \\[6pt]
 xv &  yw &  zu
 \end{array} \right| \\[6pt]
 & \hspace*{-1.05in} = \frac{1}{\Box}
 (1  +ywzu +zx -z -zu - xyw)  \\[2pt]
 & \hspace*{-0.75in} + \frac{1}{\Box}
 ( z + yw zu x  + zxy  - zx - zuy  - zxyw ) \\[2pt]
 & \hspace*{-0.5in} + \frac{1}{\Box}
 ( zu + ywzuxv  + zxyw - zxv  - zuyw - zuxyw)  \\[2pt]
 & \hspace*{-1.05in} = \frac{1}{\Box}
 (1 - x y w - x v z - u y z + x y z + u v w x y z).
\end{align*}
}
\!\!This proves \eqref{eqD}.  As a consequence of our theorem, we also obtain the following unification of the theorems of Ceva and Menelaus.

\begin{corollary*}
The points $P, Q$ and $R$ defined in \eqref{eqgKLM} are collinear if and only if
 \[
1-x y w - x v z - u y z + x y z + xyzuvw = 0.
 \]
\end{corollary*}

For an animation which ties together these two classical theorems of Euclidean geometry and which is inspired by this corollary, see \cite{Curgus}.

\section{A connection to Feynman's triangle.}

The Routh's triangle with coefficients $x=y=z=2$ attracted more attention than any other, see for example \cite[p.~9]{Steinhaus} and \cite{Cook}. The area of this special Routh's triangle equals $1/7$ of the area of its host triangle. The proof of this fact in \cite[p.~9]{Steinhaus} is a tiling type argument almost without words, while \cite{Cook} gives three different proofs using standard tools of Euclidean geometry.  According to  \cite{Cook}, this result kept Richard Feynman busy during a dinner. This story gave the triangle the name Feynman's triangle; another name for it is the ``one-seventh area triangle''.

Further on, we assume that the host triangle has area $1$. Feynman's triangle is unique in the sense that no other Routh's triangle with equal integer coefficients $x,y,$ and $z$ has an area which is a reciprocal of a nonzero digit. For example, if $x=y=z=4$, the area is $3/7$.  For all other equal integers as coefficients, the areas are more complicated fractions. In fact, even if we allow unequal positive digit coefficients $x, y$, and $z$, we can only get the area $1/2$ for $x=7,y=6,z=3$, and $1/4$ for $x = 7, y = 4,z = 1$. No other rational number in $(0,1)\setminus \{1/7, 1/4, 3/7, 1/2\}$ with digits as its  numerator and denominator is attainable as an area of Routh's triangle with positive digits as coefficients.

\begin{figure}[H]

\begin{minipage}{0.325\textwidth}
\psfrag{C}[][]{\begin{picture}(0,0)
            \put(-5,8){\makebox(0,0)[l]{$C$}}
                        \end{picture}}
\psfrag{Bv}[][]{\begin{picture}(0,0)
            \put(-16,5){\makebox(0,0)[l]{$B_u$}}
                        \end{picture}}
\psfrag{Ax}[][]{\begin{picture}(0,0)
            \put(0,6){\makebox(0,0)[l]{$A_x$}}
                        \end{picture}}
\psfrag{rA}[][]{\begin{picture}(0,0)
            \put(-6,9){\makebox(0,0)[l]{}} 
                             \end{picture}}
\psfrag{By}[][]{\begin{picture}(0,0)
            \put(-16,5){\makebox(0,0)[l]{$B_x$}}
                        \end{picture}}
\psfrag{rB}[][]{\begin{picture}(0,0)
            \put(-10,-4){\makebox(0,0)[l]{}} 
                        \end{picture}}
\psfrag{Au}[][]{\begin{picture}(0,0)
            \put(0,6){\makebox(0,0)[l]{$A_u$}}
                        \end{picture}}
\psfrag{rC}[][]{\begin{picture}(0,0)
            \put(-10,-8){\makebox(0,0)[l]{}} 
                        \end{picture}}
\psfrag{A}[][]{\begin{picture}(0,0)
            \put(-8,-8){\makebox(0,0)[l]{$A$}}
                        \end{picture}}
\psfrag{Cw}[][]{\begin{picture}(0,0)
            \put(-4,-9){\makebox(0,0)[l]{$C_u$}}
                        \end{picture}}
\psfrag{Cz}[][]{\begin{picture}(0,0)
            \put(-4,-9){\makebox(0,0)[l]{$C_x$}}
                        \end{picture}}
\psfrag{B}[][]{\begin{picture}(0,0)
            \put(-2,-8){\makebox(0,0)[l]{$B$}}
                        \end{picture}}

\setlength{\abovecaptionskip}{0pt}%
\setlength{\belowcaptionskip}{-6pt}%

\resizebox{\textwidth}{!}{\includegraphics{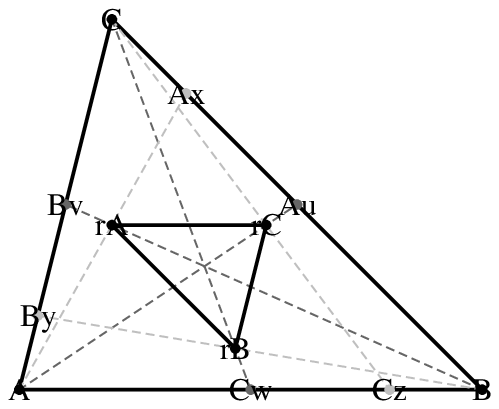}}
  \caption{} \label{f4}
\end{minipage}
\begin{minipage}{0.325\textwidth}

\psfrag{C}[][]{\begin{picture}(0,0)
            \put(-5,8){\makebox(0,0)[l]{$C$}}
                        \end{picture}}
\psfrag{Bv}[][]{\begin{picture}(0,0)
            \put(-16,5){\makebox(0,0)[l]{$B_u$}}
                        \end{picture}}
\psfrag{Ax}[][]{\begin{picture}(0,0)
            \put(0,6){\makebox(0,0)[l]{$A_x$}}
                        \end{picture}}
\psfrag{rA}[][]{\begin{picture}(0,0)
            \put(-6,9){\makebox(0,0)[l]{}} 
                             \end{picture}}
\psfrag{By}[][]{\begin{picture}(0,0)
            \put(-16,5){\makebox(0,0)[l]{$B_x$}}
                        \end{picture}}
\psfrag{rB}[][]{\begin{picture}(0,0)
            \put(-10,-4){\makebox(0,0)[l]{}} 
                        \end{picture}}
\psfrag{Au}[][]{\begin{picture}(0,0)
            \put(0,6){\makebox(0,0)[l]{$A_u$}}
                        \end{picture}}
\psfrag{rC}[][]{\begin{picture}(0,0)
            \put(-10,-8){\makebox(0,0)[l]{}} 
                        \end{picture}}
\psfrag{A}[][]{\begin{picture}(0,0)
            \put(-8,-8){\makebox(0,0)[l]{$A$}}
                        \end{picture}}
\psfrag{Cw}[][]{\begin{picture}(0,0)
            \put(-4,-9){\makebox(0,0)[l]{$C_u$}}
                        \end{picture}}
\psfrag{Cz}[][]{\begin{picture}(0,0)
            \put(-4,-9){\makebox(0,0)[l]{$C_x$}}
                        \end{picture}}
\psfrag{B}[][]{\begin{picture}(0,0)
            \put(-2,-8){\makebox(0,0)[l]{$B$}}
                        \end{picture}}

\setlength{\abovecaptionskip}{0pt}%
\setlength{\belowcaptionskip}{-6pt}%

\resizebox{\textwidth}{!}{\includegraphics{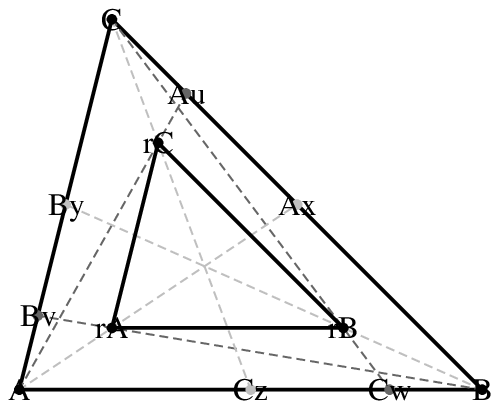}}
  \caption{} \label{f5}
\end{minipage}
\begin{minipage}{0.325\textwidth}
\psfrag{C}[][]{\begin{picture}(0,0)
            \put(-5,8){\makebox(0,0)[l]{$C$}}
                        \end{picture}}
\psfrag{Bv}[][]{\begin{picture}(0,0)
            \put(-16,5){\makebox(0,0)[l]{$B_u$}}
                        \end{picture}}
\psfrag{Ax}[][]{\begin{picture}(0,0)
            \put(0,6){\makebox(0,0)[l]{$A_x$}}
                        \end{picture}}
\psfrag{rA}[][]{\begin{picture}(0,0)
            \put(-6,9){\makebox(0,0)[l]{}} 
                             \end{picture}}
\psfrag{By}[][]{\begin{picture}(0,0)
            \put(-16,5){\makebox(0,0)[l]{$B_x$}}
                        \end{picture}}
\psfrag{rB}[][]{\begin{picture}(0,0)
            \put(-10,-4){\makebox(0,0)[l]{}} 
                        \end{picture}}
\psfrag{Au}[][]{\begin{picture}(0,0)
            \put(0,6){\makebox(0,0)[l]{$A_u$}}
                        \end{picture}}
\psfrag{rC}[][]{\begin{picture}(0,0)
            \put(-10,-8){\makebox(0,0)[l]{}} 
                        \end{picture}}
\psfrag{A}[][]{\begin{picture}(0,0)
            \put(-8,-8){\makebox(0,0)[l]{$A$}}
                        \end{picture}}
\psfrag{Cw}[][]{\begin{picture}(0,0)
            \put(-4,-9){\makebox(0,0)[l]{$C_u$}}
                        \end{picture}}
\psfrag{Cz}[][]{\begin{picture}(0,0)
            \put(-4,-9){\makebox(0,0)[l]{$C_x$}}
                        \end{picture}}
\psfrag{B}[][]{\begin{picture}(0,0)
            \put(-2,-8){\makebox(0,0)[l]{$B$}}
                        \end{picture}}

\setlength{\abovecaptionskip}{0pt}%
\setlength{\belowcaptionskip}{-6pt}%

\resizebox{\textwidth}{!}{\includegraphics{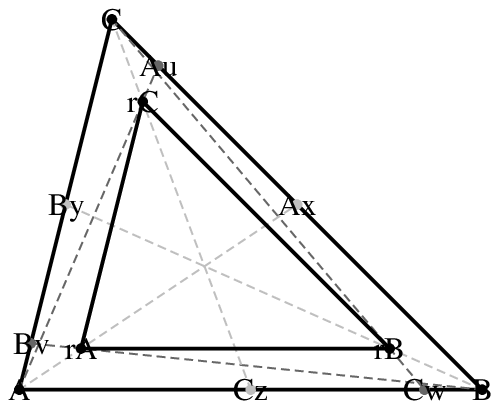}}
  \caption{} \label{f6}
\end{minipage}

\end{figure}

Does this change if we consider generalized Routh's triangles with $u=v=w$, $x=y=z$ and $x\neq u$ being digits? The answer is yes. We can get the area $1/9$ with $u=1, x=4$, see Figure~\ref{f4}, the area $1/4$ with $u=4,x=1$, see Figure~\ref{f5} and the area $4/9$ with $u=7, x =1$, see Figure~\ref{f6}. However, none of these triangles are very interesting since they are simply scaled images of the host triangle with the same center of gravity.

\begin{figure}[H]
\begin{minipage}{0.49\textwidth}

\psfrag{C}[][]{\begin{picture}(0,0)
            \put(-5,8){\makebox(0,0)[l]{$C$}}
                        \end{picture}}
\psfrag{Bv}[][]{\begin{picture}(0,0)
            \put(-16,5){\makebox(0,0)[l]{$B_u$}}
                        \end{picture}}
\psfrag{Ax}[][]{\begin{picture}(0,0)
            \put(0,6){\makebox(0,0)[l]{$A_x$}}
                        \end{picture}}
\psfrag{rA}[][]{\begin{picture}(0,0)
            \put(-6,9){\makebox(0,0)[l]{}} 
                             \end{picture}}
\psfrag{By}[][]{\begin{picture}(0,0)
            \put(-16,5){\makebox(0,0)[l]{$B_x$}}
                        \end{picture}}
\psfrag{rB}[][]{\begin{picture}(0,0)
            \put(-10,-4){\makebox(0,0)[l]{}} 
                        \end{picture}}
\psfrag{Au}[][]{\begin{picture}(0,0)
            \put(0,6){\makebox(0,0)[l]{$A_u$}}
                        \end{picture}}
\psfrag{rC}[][]{\begin{picture}(0,0)
            \put(-10,-8){\makebox(0,0)[l]{}} 
                        \end{picture}}
\psfrag{A}[][]{\begin{picture}(0,0)
            \put(-8,-8){\makebox(0,0)[l]{$A$}}
                        \end{picture}}
\psfrag{Cw}[][]{\begin{picture}(0,0)
            \put(-4,-9){\makebox(0,0)[l]{$C_u$}}
                        \end{picture}}
\psfrag{Cz}[][]{\begin{picture}(0,0)
            \put(-4,-9){\makebox(0,0)[l]{$C_x$}}
                        \end{picture}}
\psfrag{B}[][]{\begin{picture}(0,0)
            \put(-2,-8){\makebox(0,0)[l]{$B$}}
                        \end{picture}}

\setlength{\abovecaptionskip}{0pt}%
\setlength{\belowcaptionskip}{-6pt}%

\resizebox{\textwidth}{!}{\includegraphics{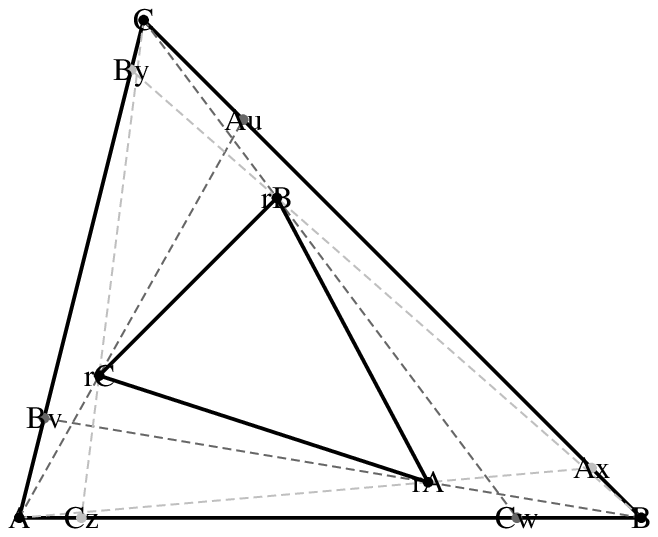}}
  \caption{$1/4$ area; $u=4,x=1/9$} \label{f7}
\end{minipage}
\begin{minipage}{0.49\textwidth}
\psfrag{C}[][]{\begin{picture}(0,0)
            \put(-5,8){\makebox(0,0)[l]{$C$}}
                        \end{picture}}
\psfrag{Bv}[][]{\begin{picture}(0,0)
            \put(-16,5){\makebox(0,0)[l]{$B_u$}}
                        \end{picture}}
\psfrag{Ax}[][]{\begin{picture}(0,0)
            \put(0,6){\makebox(0,0)[l]{$A_x$}}
                        \end{picture}}
\psfrag{rA}[][]{\begin{picture}(0,0)
            \put(-6,9){\makebox(0,0)[l]{}} 
                             \end{picture}}
\psfrag{By}[][]{\begin{picture}(0,0)
            \put(-16,5){\makebox(0,0)[l]{$B_x$}}
                        \end{picture}}
\psfrag{rB}[][]{\begin{picture}(0,0)
            \put(-10,-4){\makebox(0,0)[l]{}} 
                        \end{picture}}
\psfrag{Au}[][]{\begin{picture}(0,0)
            \put(0,6){\makebox(0,0)[l]{$A_u$}}
                        \end{picture}}
\psfrag{rC}[][]{\begin{picture}(0,0)
            \put(-10,-8){\makebox(0,0)[l]{}} 
                        \end{picture}}
\psfrag{A}[][]{\begin{picture}(0,0)
            \put(-8,-8){\makebox(0,0)[l]{$A$}}
                        \end{picture}}
\psfrag{Cw}[][]{\begin{picture}(0,0)
            \put(-4,-9){\makebox(0,0)[l]{$C_u$}}
                        \end{picture}}
\psfrag{Cz}[][]{\begin{picture}(0,0)
            \put(-4,-9){\makebox(0,0)[l]{$C_x$}}
                        \end{picture}}
\psfrag{B}[][]{\begin{picture}(0,0)
            \put(-2,-8){\makebox(0,0)[l]{$B$}}
                        \end{picture}}

\setlength{\abovecaptionskip}{0pt}%
\setlength{\belowcaptionskip}{-6pt}%

\resizebox{\textwidth}{!}{\includegraphics{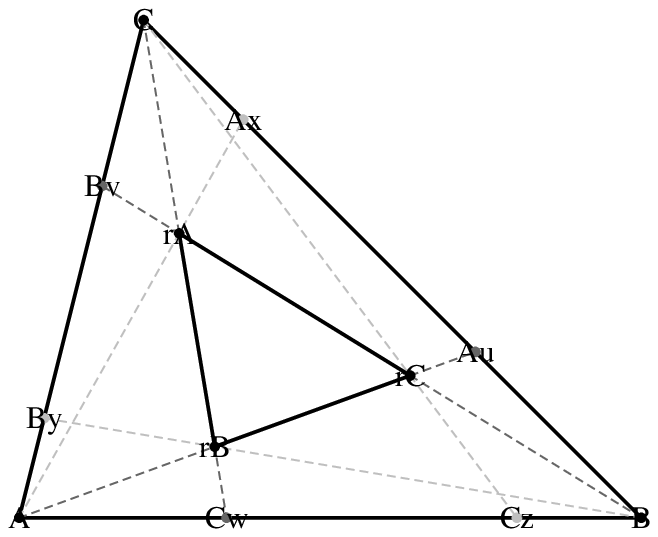}}
  \caption{$1/7$ area; $u=1/2,x=4$} \label{f8}
\end{minipage}
\end{figure}

Nevertheless, if we allow one of the coefficients $u$ or $x$ to be a reciprocal of a nonzero digit, then, besides the previous three configurations, we can also get the area $1/4$ with $u=9, x =1/4$ or $u=4,x=1/9$, and our old acquaintance, the area $1/7$, with $u = 4,x=1/2$, or $u=2,x=1/4$, or $u=1/2,x=4$; see Figures~\ref{f7} and~\ref{f8}.

\vspace*{-8pt}

We conclude this note with an invitation for the reader to find a purely geometric argument along the lines of \cite[p.~9]{Steinhaus} for the claims indicated in the last two figures.

\bigskip

\noindent\textit{Department of Mathematics, Western Washington University, Bellingham, WA 98225, USA, \\ Arpad.Benyi@wwu.edu}

\bigskip

\noindent\textit{Department of Mathematics, Western Washington University, Bellingham, WA 98225, USA, \\ Branko.Curgus@wwu.edu}

\end{document}